\documentclass[review,onefignum,onetabnum]{siamart190516}

\usepackage{amsmath}
\usepackage{listings}
\usepackage{booktabs}

\makeatletter
\def\myunderbrace#1{\mathop{\vtop{\m@th\ialign{##\crcr
   $\hfil\displaystyle{#1}\hfil$\crcr
   \noalign{\kern3\p@\nointerlineskip}%
   \footnotesize\downbracefill\crcr\noalign{\kern3\p@\nointerlineskip}}}}\limits}
\makeatother

\graphicspath{{../figures/}}

\usepackage{lipsum}
\usepackage{amsfonts}
\usepackage{mathrsfs}
\usepackage{graphicx}
\usepackage{epstopdf}
\usepackage{xcolor} 
\usepackage{todonotes}
\usepackage{algorithm}
\usepackage{algpseudocode}
\usepackage[caption=false]{subfig}
\usepackage{multirow}
\usepackage{MnSymbol}
\graphicspath{{figures/}}
\ifpdf
  \DeclareGraphicsExtensions{.eps,.pdf,.png,.jpg}
\else
  \DeclareGraphicsExtensions{.eps}
\fi
\usepackage{float} 

\usepackage{cleveref}

\newsiamremark{remark}{Remark}
\newsiamremark{hypothesis}{Hypothesis}
\newsiamremark{example}{Example}
\crefname{hypothesis}{Hypothesis}{Hypotheses}
\newsiamthm{claim}{Claim}

\headers{A new banded Petrov--Galerkin spectral method}{Ouyuan Qin, Lu Cheng, and Kuan Xu}

\title{A new banded Petrov--Galerkin spectral method}

\author{Ouyuan Qin\thanks{School of Mathematical Sciences, University of Science and Technology of China, 96 Jinzhai Road, Hefei, Anhui, 230026, China (\email{ouyuanqin@mail.ustc.edu.cn}, \email{chengllu@mail.ustc.edu.cn}, \email{kuanxu@ustc.edu.cn}).}
\and Lu Cheng\footnotemark[1]
\and Kuan Xu\footnotemark[1]}
\usepackage{amsopn}
\DeclareMathOperator{\diag}{diag}

\newcommand{\bsu}{\boldsymbol{u}}

\newcommand{\ha}{\hat{a}}
\newcommand{\hb}{\hat{b}}
\newcommand{\hc}{\hat{c}}

\newcommand{\mB}{\mathcal{B}}
\newcommand{\mD}{\mathcal{D}}

\newcommand{\mI}{\mathcal{I}}

\newcommand{\mL}{\mathcal{L}}
\newcommand{\mM}{\mathcal{M}}
\newcommand{\mO}{\mathcal{O}}
\newcommand{\mP}{\mathcal{P}}

\newcommand{\mS}{\mathcal{S}}

\def\bf{\boldsymbol{f}}

\def\diag{\mathop{\mathrm{diag}}}

\def\ha{\hat{a}}

\begin{document}

\maketitle

\begin{abstract}
We propose a Petrov--Galerkin spectral method for ODEs with variable coefficients. When the variable coefficients are smooth, the new method yields a strictly banded linear system, which can be efficiently constructed and solved in linear complexity. The performance advantage of our method is demonstrated through benchmarking against Mortensen's Galerkin method and the ultraspherical spectral method. Furthermore, we introduce a systematic approach for designing the recombined basis and establish that our new method serves as a unifying framework that encompasses all existing banded Galerkin spectral methods. This significantly addresses the ongoing challenge of developing recombined bases and sparse Galerkin spectral method. Additionally, the accelerating techniques presented in this paper can also enhance the performance of the ultraspherical spectral method.
\end{abstract}
  
  

\section{Introduction}\label{sec:intro}
In this paper, we propose a general framework for constructing fast Petrov--Galerkin (PG) spectral methods that solve the ordinary differential equation
\begin{subequations}
\begin{align}
\mathcal{L} u(x) &= g(x),~~~ x \in [-1 , 1], \label{lug}\\
\text{s.t. } \mathcal{B} u(x) &= c, \label{buc}
\end{align} \label{intode}%
\end{subequations}
where $\mathcal{B} = (\mathcal{B}_0, \mathcal{B}_1, \dots, \mathcal{B}_{N-1})^{\top}$ contains $N$ linear constraints, including boundary conditions or side constraints or a mix of them, and $c$ is an $N$-vector. We assume that the differential operator
\begin{align}
\mathcal{L} = a^N(x)\frac{{\mathrm{d}}^N}{{\mathrm{d}}x^N} + \ldots + a^1(x)\frac{{\mathrm{d}}}{{\mathrm{d}}x} + a^0(x) \label{mL}
\end{align} 
for $a^{N}(x) \neq 0$ and the variable coefficients $a^0(x), \ldots, a^N(x)$ and the right-hand side function $g(x)$ possess certain regularity. 

There has been a longstanding effort to develop basis and weight functions that ensure Galerkin spectral methods yield sparse or structured linear systems, enabling efficient solutions. In an earliest attempt \cite{hei1}, a set of recombined Chebyshev polynomials are used as trial functions for significant improvement in the conditioning of the discrete systems. In passing, it is realized that the underlying rank structure of the resulting systems allows for fast solution. The seminal paper \cite{she1} by Shen proposed a recombined Legendre basis for a Bubnov--Galerkin (BG) spectral method which leads to banded matrices for second- and fourth-order ODEs with constant coefficients, while the Chebyshev version was explored in \cite{she2} yielding low-rank upper Hessenberg systems which can also be solved in a linear complexity. Doha extended Shen's methods to ultraspherical and general Jacobi polynomials, addressing both odd- and even-order ODEs in a series of papers \cite{doh3,doh4,doh2,doh1}. All these works, however, concentrate on ODEs with constant coefficients. The first attempt in this vein towards variable-coefficient ODEs is Mortensen's Petrov--Galerkin (MPG) spectral method \cite{mor1}, which also gives rise to banded systems. But it is only for variable coefficients in the form of power series centered at the origin, i.e., $\sum_{k} c_k x^k$, the linear system of MPG method can be constructed in a complexity that is proportional to the size of discretization. For ODEs with coefficients of general univariate functions, MPG method has to resort to numerical integration which results in a quadratic complexity. Another relevant work \cite{she3} by Shen, Wang, and Xia, instead of focusing on constructing a banded Galerkin method, shows that the linear systems arising from BG spectral methods, although usually full, can be solved with reduced complexity by exploring the underlying rank structure. \cref{tab:cost_pgm_cont} summarizes the complexities of these methods in terms of constructing the coefficient matrix and obtaining the solution.

\begin{table}[t!]
  \centering
  \renewcommand{\arraystretch}{1.4}
  \setlength{\tabcolsep}{8pt}
  \caption{Cost of Galerkin spectral methods by construction and solution. $N$ and $n$ are the order of the ODE and the dimension of the discretized system respectively. For the first 4 rows, the complexity is for a second-order ODE. The cost of constructing MPG is $\mO(n)$ using recursion and $\mO(n^2)$ using numerical integration.}
  \begin{tabular}{l lcc}
  \toprule
  coefficients & method & construction & solution \\
  \midrule
  \multirow{6}{*}{constant} 
  & Heinrichs \cite{hei1} & $\frac{1}{2}n^2 + \frac{1}{2}n + 2$ & $11n - 20$ \\
  & Shen \cite{she1} & $9n-12$ & $7n-12$ \\
  & Shen \cite{she2} & $\frac{1}{2}n^2 + \frac{21}{2}n - 5$ & $11n - 20$ \\
  & Elbarbary \cite{elb1} & $50n - 58$ &  $11n - 12$ \\
  & Doha and Abd-Elhameed \cite{doh5} & $\mO(Nn^2)$ & $\frac{2}{3}(N+9)n^2$ \\
  & Doha \cite{doh4} & $\mO(N^2n)$ & $\frac{2}{3}(N+7)Nn$ \\
  \midrule
  \multirow{2}{*}{variable} 
  & Shen, Wang and Xia \cite{she3} & $\mO(n\log^2 n)$ & $\mO(n)$ \\
  & Mortensen \cite{mor1} & $\mO(n)$ or $\mO(n^2)$ & $\mO(n)$ \\
  \bottomrule
  \end{tabular}
  \label{tab:cost_pgm_cont}
\end{table}

The new PG spectral method that we propose distinguishes itself from existing methods in a few perspectives: (1) With recombined basis and weight functions that are carefully designed, it leads to a strictly banded system for general variable coefficients, provided that the variable coefficients can be approximated by series of classical orthogonal polynomials, e.g., Chebyshev series; (2) For such a banded system, both the construction and the solution costs are linearly proportional to the size $n$ of the system. Particularly, the banded systems can exploit the standard library subroutines to gain more advantage in speed; (3) It can be shown to serve as an overarching method for all existing banded Galerkin methods. 

The Galerkin methods are sometimes argued to be difficult to automate for boundary conditions, and the design of the recombined basis is deemed as more of an art than a science. To alleviate the pain of basis design, we propose a systematic approach to recombining basis. When this approach is implemented symbolically, the stencil coefficients it produces coincide, up to a scaling factor, with those recommended in the aforementioned studies.

Meanwhile, the techniques introduced for the new PG method can also be applied to the ultraspherical spectral (US) method to have it significantly accelerated.

Throughout this paper, we shall make frequent use of quasimatrices. For $x\in [a, b]$, $\left(g_1(x)| g_2(x) | \cdots| g_n(x)\right)$ is an $[a, b] \times n$ column quasi-matrix, which is in fact a matrix with $n$ columns. Each column is a univariate function defined on an interval $[a, b]$, and can be deemed as a continuous analogue of a tall-skinny matrix, where the rows are indexed by a continuous, rather than discrete, variable. A row quasi-matrix is the transpose of a column one. In this paper, the unitary and bilinear operations of quasimatrices follow exactly those of standard vectors, including scalar multiplication and outer product. For the notion of quasi-matrices, see, for example, \cite{bat1}.

To facilitate the exposition, we may omit the argument $x$ of a function or a polynomial when it is clear from the context. For example, we may write a function $g(x)$ as simple as $g$.

The article is organized as follows. In \cref{sec:bpg}, we present the method for designing recombined basis functions that satisfy general linear constraints. With the recombined Chebyshev and ultraspherical polynomials the proposed PG method is shown to produce banded linear systems. In \cref{sec:fs}, we show how such a banded system can be constructed in a linear complexity. The advantage in speed is demonstrated in \cref{sec:exm}, where the proposed method is compared with MPG \cite{mor1} and the US \cite{olv1} methods. In \cref{sec:jacobi}, the proposed method is extended to Jacobi polynomials to show that MPG method and other earlier sparse Galerkin methods are specific instances of the proposed framework. \Cref{sec:us} demonstrates that the techniques introduced in \cref{sec:bpg,sec:fs} can also accelerate the US method. We close by a discussion.


\section{A banded PG method} \label{sec:bpg}
Throughout this paper, we assume that the linear constraints \cref{buc} are homogeneous, i.e., $c=0$. Otherwise, the lifting technique \cite[\S 8.1]{she4} can be used to convert the problem to a homogeneous one by subtracting from the solution a low degree polynomial that satisfies the inhomogeneous linear constraints. 

Let $P_n$ be the space consisting of polynomials of degree less than or equal to $n$. We denote the trial space by $V_n = \{ v \in P_{n + N - 1} : \mathcal{B}v = 0 \}$ and the test space by $W_n$. $W_n$ is usually of dimension $n$ and a subspace of $P_{n + N - 1}$ whose elements satisfy certain given constraints, such as boundary conditions. For odd $ N $, a preferable choice is $ W_n = V^*_n $, where the elements in $V^*_n$ satisfy the dual boundary condition \cite{she5}. For even $N$, a common practice is $W_n = V_n$ \cite{she1, she2}. Given a weight function $\omega(x)$, the standard PG method seeks a solution $u_n \in V_n$ so that
\begin{align}
(\mathcal{L} u_n(x) - g(x), w(x))_{\omega} = 0,~~~~~\forall w(x) \in W_n, \label{pgsm}
\end{align}
where $(\cdot, \cdot)_{\omega}$ denotes an inner product with respect to $\omega$. Let 
\begin{align*}
V_n &= \{ \phi_0(x), \phi_1(x), \ldots, \phi_{n-1}(x) \}, \\
W_n &= \{ \psi_0(x), \psi_1(x), \ldots, \psi_{n-1}(x) \},
\end{align*}%
where $\phi_k(x)$ and $\psi_k(x)$ are the trial and test functions respectively. If we assume that $u_n = \sum_{k=0}^{n-1} v_k \phi_k$, \cref{pgsm} leads to
\begin{align}
A v = f,  \label{pgsm_mtx}
\end{align}
where 
\begin{align*}
v = (v_0, v_1, \ldots, v_{n-1}), ~~~f = (f_0, f_1, \ldots, f_{n-1}), ~~~ f_i = (\psi_i, g)_{\omega},
\end{align*}
and the $(i,j)$th entry of $A$
\begin{align}
A(i,j) = (\psi_i, \mathcal{L} \phi_j)_{\omega}. \label{Aij}
\end{align}
Often the entries of $A$ can only be obtained by numerical integration, except in certain simplest scenarios of constant-coefficient ODEs where the entries can be spelled out in closed forms. We shall nonetheless show that the construction of $A$ in the proposed framework involves neither quadrature nor manual calculation.

\subsection{Trial and test bases}\label{sec:tt}
The framework we present in this paper can be constructed using Jacobi polynomials. However, we begin our discussion with Chebyshev and ultraspherical polynomials and defer the generalization to Jacobi to \cref{sec:jacobi}. In principle, there are infinitely many ways to combine Chebyshev polynomials so that $\phi_k(x)$ satisfies the linear constraints. To keep the resulting system as banded as possible, it is however preferable to combine only a small number, say, e.g., $N+1$, of the consecutively neighboring Chebyshev polynomials. That is,
\begin{align}
  V_n = \left(\phi_0(x)\,\middle|\, \phi_1(x)\,\middle|\, \cdots\,\middle|\, \phi_{n-1}(x) \right) = \left(T_0(x)\,\middle|\, T_1(x)\,\middle|\, \cdots\,\middle|\, T_{n+N-1}(x) \right) R, \label{baseT}
\end{align}
where the stencil matrix
\begin{align}
R = 
  \begin{pmatrix}
    \gamma_{0}^{0} &           &      &      \\
    \gamma_{1}^{0} & \gamma_{0}^{1} &      &      \\
    \vdots      & \gamma_{1}^{1} & \ddots &      \\
    \gamma_{N}^{0} & \vdots      & \ddots & \gamma_{0}^{n-1} \\
                & \gamma_{N}^{1} & \ddots & \gamma_{1}^{n-1} \\
                &             & \ddots & \vdots \\
                &             &        & \gamma_{N}^{n-1}
  \end{pmatrix} 
  \in \mathbb{R}^{(n + N - 1) \times n}. \label{stencil}
\end{align}
Let $B_{ij}^{k} = \mathcal{B}_{i}(T_{k+j}(x))$, where $\mathcal{B}_i$ is the $i$th linear constraint. For the $k$th column of $R$, the recombination coefficients $\{\gamma_{j}^{k}\}_{j=0}^N$ can be determined by solving the under-determined system
\begin{align}
  \label{coe_bc}
  \begin{pmatrix}
    B_{00}^{k} & B_{01}^{k} & \ldots & B_{0N}^{k} \\
    B_{10}^{k} & B_{11}^{k} & \ldots & B_{1N}^{k} \\
    \vdots     & \vdots     & \vdots & \vdots     \\
    B_{N-1,0}^{k} & B_{N-1,1}^{k} & \ldots & B_{N-1,N}^{k}   
  \end{pmatrix}
  \begin{pmatrix}
    \gamma_{0}^{k}\\
    \gamma_{1}^{k}\\
    \vdots     \\
    \gamma_{N}^{k}
  \end{pmatrix} = 
  \begin{pmatrix}
    0 \\
    0 \\
    \vdots \\
    0
  \end{pmatrix},
\end{align}
which amounts to requiring the $k$th recombined basis function to satisfy the homogeneous linear constraints. In \cref{coe_bc}, the total number of the unknowns is $N+1$, greater than that of the equations by one. For a stencil that combines only $N$ Chebyshev polynomials, the resulting system  would only have a unique but trivial solution $\gamma_{j}^{k}=0$. Hence, $N+1$ is the smallest number to guarantee a nontrivial recombination. In additional, it also leads to a minimal upper bandwidth as we shall see in \cref{sec:bcm}. Mathematically, $\{\gamma_{j}^{k}\}_{j=0}^N$ can be any value that satisfies \cref{coe_bc} as long as the highest degree coefficient $\gamma_{N}^{k} \neq 0$. Numerically, $\gamma_{N}^{k}$ should not be too small for $V_n$ to be a stable basis. In practice, the extra degree of freedom is exploited by setting $\gamma_{N}^{k}$ to a value bounded below by a safeguard value $\gamma_{\min}$, leaving the other $\{\gamma_{j}^{k}\}_{j=0}^{N-1}$ to be determined uniquely. 

\begin{algorithm}[!t]
  \caption{Algorithm for calculating the analytic expression for $\gamma_j^k$}
  \label{alg:bcs}
  \begin{description}
      \item[Inputs:] The total number $N$ of linear constraints and the functionals $\mathcal{B}_{i}$, $i = 0, \dots, N-1$ that maps a Chebyshev polynomial to a number.
      \item[Outputs:] The analytic expression for $\gamma_j^k$, $j = 0, \dots, N$.
  \end{description}
  \vspace{2.5pt}
  \hrule
  \begin{algorithmic}[1]
      \State Initialize an empty set $E$ for equations.
      \For{$i = 0, \dots, N-1$} \Comment{Loop over the rows}
        \State $e_i = 0$ 
        \For{$j = 0, \dots, N$} \Comment{Loop to set up each row of \eqref{coe_bc}}
          \State Find the expression $B = \mathcal{B}_{i}(T_{k+j})$
          \State $e_i = e_i + B \gamma_{j}^k$
        \EndFor
        \State Append the $i$th equation $e_i = 0$ to set $E$.
      \EndFor
      \State Replace $\gamma_{N}^k$ by $1$. \Comment{$\gamma_{N}^k$ is anchored for the solution to be unique.}
      \State Solve the system $E$ for $\gamma_{j}^k$, $j = 0, \dots, N-1$.
  \end{algorithmic}
\end{algorithm}

Since the explicit expressions for $B_{ij}^k$ are usually available for Chebyshev polynomials, we recommend \cref{coe_bc} be populated with these expressions and solved symbolically to produce closed-form formulae for $\{\gamma_{j}^{k}\}_{j=0}^{N-1}$. Note that \cref{coe_bc} is solved only once and the solution is symbolic expressions for $\{\gamma_{j}^{k}\}_{j=0}^{N-1}$ in terms of $k$. The procedure is summarized in \cref{alg:bcs}. We also include a short \textsc{Mathematica} program in \cref{sec:mcode} for an example of \cref{buc} with the linear constraints $\mathcal{B}_{i}$ all being endpoint boundary conditions of proper orders. The expressions for the combination coefficients $\gamma_j^k$ obtained by this \textsc{Mathematica} program are exactly those used in \cite{jul,mor}, up to a scaling factor.

In principal, one can also evaluate $B_{ij}^k$ and solve \cref{coe_bc} numerically. This is however impractical for two reasons. First, unlike symbolic computation which works with expressions so that we solve \cref{coe_bc} only once, we have to set up and solve \cref{coe_bc} $n$ times for $k = 0, 1, \ldots, n-1$ if we switch to numerical computation. This is not cheap. Second and more importantly, since in many cases it is almost inevitable that \cref{coe_bc} is ill-conditioned, the numerical solution of $\{\gamma_{j}^{k}\}_{j=0}^{N-1}$ could be very inaccurate, if not totally erroneous. For example, the boundary values of the $p$th derivative of Chebyshev polynomial $T_k(x)$ are $\mO(k^{2p})$. Thus, if some of the linear constraints in \cref{buc} are high-order boundary conditions, the poor conditioning of \cref{coe_bc} would prevent the numerical computation from producing any meaningful solutions.

As we shall see shortly below, the test basis functions are chosen to be the combinations of a small number of neighboring ultraspherical polynomials for the method to be banded. That is,
\begin{align}
\Psi = \left(\psi_0(x) \,\middle|\, \psi_1(x) \,\middle|\, \cdots \,\middle|\, \psi_{n-1}(x) \right) = \left(C^{(N)}_0(x) \,\middle|\, C^{(N)}_1(x) \,\middle|\, \cdots \,\middle|\, C^{(N)}_{n+N-1}(x) \right) Q. \label{baseCN}
\end{align}  
Like for the trial functions, we determine the stencil matrix $Q \in \mathbb{R}^{(n+N-1) \times n}$ also by \cref{alg:bcs} but with Chebyshev polynomial $T_{k+j}$ replaced by $C^{(N)}_{k+j}$. This way, the recombined test functions also satisfy the homogeneous linear constraints and keep the lower bandwidth as small as possible. If the test functions are chosen to satisfy the dual boundary conditions, the linear constraint $\mB_i$ in line 5 of \cref{alg:bcs} should be changed accordingly. 

A more advanced \textsc{Mathematica} program \texttt{BasisRecombination.nb} capable of handling a broader range of constraints, e.g., midpoint conditions and global condition of integration, is available online from \cite{qin2}.

\subsection{Banded coefficient matrix}\label{sec:bcm}
Our discussion in the rest of this section makes use of the differentiation, multiplication, and the conversion operators introduced in \cite{olv1} for the US method. For $k \geq 1$, the $k$th-order differentiation operator of infinite dimensions is given by
\begin{align*}
\begin{aligned}
  & ~~\,\,\, \myunderbrace{\text{\footnotesize $k$ times}} \\[-7pt]
  \mD_{k} = 2^{k - 1}(k - 1)! & \begin{pmatrix}
       \  0 ~~  \cdots ~~ 0  &k & & & \\
      & &k +1& & \\
      & & &k +2 & \\
      & & & &\ddots\\
  \end{pmatrix}
\end{aligned},
\label{D}
\end{align*}
which maps an infinite vector of Chebyshev $T$ coefficients to the coefficient vector of its $k$th derivative but in $C^{(k)}$. If $a^0(x)$  is written as an infinite Chebyshev series, i.e., $a^0(x)= \sum_{j=0}^{\infty}a^0_j T_j(x)$, the action of multiplying another Chebyshev series by $a^0(x)$ is represented by an almost Toeplitz-plus-Hankel multiplication operator
\begin{align}
  \mM_{0}[a^0]=\frac{1}{2}\left[\begin{pmatrix}
  2 a^0_{0} & a^0_{1} & a^0_{2} & a^0_{3} & \cdots \\
  a^0_{1} & 2 a^0_{0} & a^0_{1} & a^0_{2} & \ddots \\
  a^0_{2} & a^0_{1} & 2 a^0_{0} & a^0_{1} & \ddots \\
  a^0_{3} & a^0_{2} & a^0_{1} & 2 a^0_{0} & \ddots \\
  \vdots & \ddots & \ddots & \ddots & \ddots
  \end{pmatrix}+\begin{pmatrix}
  0 & 0 & 0 & 0 & \cdots \\
  a^0_{1} & a^0_{2} & a^0_{3} & a^0_{4} & \cdots \\
  a^0_{2} & a^0_{3} & a^0_{4} & a^0_{5} & \udots \\
  a^0_{3} & a^0_{4} & a^0_{5} & a^0_{6} & \udots \\
  \vdots & \udots & \udots & \udots & \udots
  \end{pmatrix}\right].
  \label{M0}
\end{align}
Since the range of $\mD_{k}$ is the ultraspherical space of $C^{(k)}(x)$ for $k \geq 1$, it is natural to assume that the variable coefficients $a^k(x)$ are approximated by infinite $C^{(k)}$ series, i.e., $a^k(x) = \sum_{j=0}^{\infty} a^k_j C^{(k)}_j(x)$. Multiplying by $a^k(x)$ can then be effected by the multiplication operator $\mM_{k}[a^k]$, which maps between two $C^{(k)}$ series. The most straightforward way to generate $\mM_{k}[a^k]$ \cite[\S 6.3.1]{tow1} is to express it by a series
\begin{align*}
\mM_{k}[a^k]=\sum_{j=0}^{\infty} a^k_{j} \mM_{k}[C_{j}^{(k)}], \label{Ml}
\end{align*}
where $\mM_{k}[C_{j}^{(k)}]$ is obtained by the three-term recurrence relation
\begin{subequations}
\begin{align}
\mM_{k}[C_{j+1}^{(k)}]=\frac{2(j+k)}{j+1} \mM_{k}[x] \mM_{k}[C_{j}^{(k)}]-\frac{j+2 k-1}{j+1} \mM_{k}[C_{j-1}^{(k)}], \quad j \geq 1
\end{align}
with $\mM_{k}[C_{0}^{(k)}]$ is infinite identity operator $\mI$, $\mM_{k}[C_{1}^{(k)}]=2 k \mM_{k}[x]$, and
\begin{align}
\mM_{k}[x]=\left(\begin{array}{ccccc}
0 & \frac{2 k}{2(k+1)} & & & \\
\frac{1}{2 k} & 0 & \frac{2 k+1}{2(k+2)} & & \\
& \frac{2}{2(k+1)} & 0 & \frac{2 k+2}{2(k+3)} & \\
& & \frac{3}{2(k+2)} & 0 & \ddots \\
& & & \ddots & \ddots
\end{array}\right).
\end{align}\label{M}%
\end{subequations}
Another way to construct these multiplication operators is by an explicit but intricate formula given by equation (3.6) in \cite{olv1}. Since $a^0(x), \ldots, a^N(x)$ are assumed to be smooth, they can be approximated by finite Chebyshev or ultraspherical series to machine precision. Suppose that the approximant to $a^k (x)$ is of degree $m_{k}$. As long as $n \gg m_{k}$, $\mM_{k}[a^k]$ is a banded matrix.

The conversion operators are employed to upgrade an ultraspherical space of low order to higher ones. The transform from Chebyshev to $C^{(1)}(x)$ is effected by
\begin{align*}
\mS_0 = \begin{pmatrix}
1&  &-\frac{1}{2}&  & & \\
 &\frac{1}{2}& &-\frac{1}{2}& & \\
 & &\frac{1}{2}& &-\frac{1}{2} & \\
 & & &\ddots& &\ddots\\
\end{pmatrix}, \nonumber
\end{align*}
whereas that from $C^{(k)}(x)$ to $C^{(k+1)}(x)$ is done by
\begin{align*}
\mS_{k} = \begin{pmatrix}
1&  &-\frac{k}{k +2}&  & &\\
 &\frac{k}{k +1}& &-\frac{k}{k +3}& & \\
 & &\frac{k}{k +2}& &-\frac{k}{k +4} & \\
 & & &\ddots& &\ddots\\
\end{pmatrix}
~\text{for }k \geq 1.  \nonumber
\end{align*}

With these operators, $\mL$ given by \cref{mL} can be represented as
\begin{align}
\mathcal{M}_N[a^N] \mathcal{D}_N + \sum_{k = 0}^{N-1}{\mathcal{S}}_{N-1}\ldots {\mathcal{S}}_{k}\mathcal{M}_{k}[a^{k}]\mathcal{D}_{k},  \label{mL2}
\end{align}
where $\mD_0=\mI$. This is a well-known result from the US method. Now the first main result of this paper is in order.
\begin{theorem} \label{lem:banded_cmtx}
For the PG method with trial functions $\phi_{k}(x)$, test functions $\psi_{k}(x)$, and the weight $\omega(x) = (1 - x^2)^{N - 1 / 2}$, the coefficient matrix
\begin{align}
A = Q^{T} \Omega_{n + N - 1} L_{n + N - 1} R, \label{coe_banded_pg} 
\end{align}
where
\begin{align*}
  \Omega_{n + N - 1} = \diag 
  \left(d_0, \ldots, d_{n+N-1} \right), \; d_{j} = \frac{\pi 2^{1 - 2N}\Gamma(j + 2N)}{j! (j + N)[\Gamma(N)]^2}
\end{align*}
and 
\begin{align}
L_{n + N - 1} = \mathcal{P}_{n+N-1}\left( \mathcal{M}_N[a^N] \mathcal{D}_N + \sum_{k = 0}^{N-1}{\mathcal{S}}_{N-1}\ldots {\mathcal{S}}_{k}\mathcal{M}_{k}[a^{k}]\mathcal{D}_{k} \right) \mathcal{P}_{n+N-1}^{\top}. \label{opL_f}
\end{align}
Here, the projection operator $\mathcal{P}_{n+N-1} = \left(I_{n+N-1}, 0_{(n+N-1) \times \infty} \right)$.
\end{theorem}
\begin{proof}
First, note that
\begin{align*}
  \mathcal{L}\left(T_0(x) \,\middle|\, T_1(x) \,\middle|\, \cdots \,\middle|\, T_{n + N - 1}(x)\right) = \left(C^{(N)}_0(x) \,\middle|\, C^{(N)}_1(x) \,\middle|\, \cdots \,\middle|\, C^{(N)}_{n + N - 1}(x)\right) L_{n + N - 1}.
\end{align*}  
Substituting \cref{baseT}, \cref{baseCN}, and the last equation into \cref{Aij} gives
\begin{align*}
  A = Q^{T}
  \left(
  \begin{pmatrix}
  C^{(N)}_0(x)    \\
  C^{(N)}_1(x)    \\
  \vdots \\
  C^{(N)}_{n + N -1}(x)
  \end{pmatrix}
  \left(C^{(N)}_0(x) \,\middle|\, C^{(N)}_1(x) \,\middle|\, \cdots \,\middle|\, C^{(N)}_{n + N - 1}(x)\right) \right)_{\omega} L_{n + N - 1} R.
\end{align*}  
The orthogonality of the ultraspherical polynomials implies \cref{coe_banded_pg}.
\end{proof}

As shown above, $Q$ and $R$ are both banded and $\Omega_{n + N - 1}$ is diagonal. Since all the operators in \cref{mL2} are banded, so is $L_{n + N - 1}$ \cite{olv1}. If $m = \max\limits_{k = 0, 1, \ldots, N} \left\{m_{k}\right\}$, \cref{lem:banded_cmtx} states that for $n \gg m$ the coefficient matrix $A$ in \cref{pgsm_mtx} is strictly banded as a consequence of its being the sum of the products of a series of banded matrices.

Both LU and QR methods can be used to solve a banded system in a linear complexity. To unleash the full potential of the banded systems, it is preferable to solve \cref{pgsm_mtx} by calling \textsc{Lapack}'s subroutines \texttt{gbtrf} and \texttt{gbtrs}. In \textsc{Julia}, one also can call \texttt{qr}, which is, as far as we see, equally fast as \texttt{gbtrf} and \texttt{gbtrs}. As shown in \cref{sec:exm,sec:us}, solving such a banded system by these standard library subroutines could be much faster than solving an almost-banded system using user-supplied code, which suggests basis recombination is crucial to the performance. 



\section{Fast construction of the linear system} \label{sec:fs}
We have demonstrated that \cref{pgsm_mtx} is banded and its solution only costs $O(n)$ flops, and now turn to the construction of \cref{pgsm_mtx}. Since the Chebyshev coefficients of the variable coefficients $a^k(x)$ and the right-hand side $f$ can be obtained efficiently via FFT, we assume that they are available. In this section, we show how $A$ in \cref{pgsm_mtx} can be constructed in an $O(N^2mn)$ complexity. 

\cref{lem:banded_cmtx} shows that the construction of $A$ amounts to those of $Q$, $R$, $\Omega_{n+N-1}$, and $L_{n+N-1}$ separately before calculating the product. With the expressions for the stencil parameters $\{\gamma_{j}^{k}\}_{j=0}^N$, $Q$ and $R$ can be constructed in $\mO(Nn)$ flops. The construction of $\Omega_{n+N-1}$ also incurs $O(n)$ flops. The problem now boils down to the construction of $L_{n+N-1}$. We shall show below that an approximation of $L_{n+N-1}$ can be constructed in $O(N^2mn)$ flops. On the face of it, this is no wonder, as it is a cost proportional to $n$. This cost can be readily deduced from \cite{olv1}, even though it is not explicitly stated. However, our focus here is on the linear dependence of $m$, as we shall see below.

In practice, unlike \cref{opL_f} the truncation of \cref{mL2} is done operatorwise at each $\mD_{k}$, $\mS_{k}$, and $\mM_{k}$. Instead of exact truncation \cite[Remark 2]{olv1}, we take only square truncations of them for simplicity. For the detail of truncating \cref{mL2} exactly, see \cite[\S 3.3]{qin1}. For notational convenience, let $D_{k} = \mP_{n+N-1}\mD_{k}\mP_{n+N-1}^{\top}$, $S_{k} = \mP_{n+N-1}\mS_{k}\mP_{n+N-1}^{\top}$, and $M_{k}[a^{k}] = \mP_{n+N-1}\mM_{k}[a^{k}]\mP_{n+N-1}^{\top}$. We then approximate $L_{n+N-1}$ by
\begin{align}
L = M_N[a^N] D_N + \sum_{k = 0}^{N-1} S_{N-1} \ldots S_{k}M_{k}[a^{k}] D_{k}. \label{L}
\end{align}

The cost of constructing $D_{k}$ and $S_{k}$ is $\mO(n)$ flops and therefore minimal, since their entries are known explicitly. What remains is the construction of $M_{k}$ for $k = 0, 1, \ldots, N$. Note that constructing $M_{k}$ using either the recursive method \cite[\S 6.3.1]{tow1} or the explicit formula \cite{olv1} costs $O(m^2n)$ flops. Noting that it is usually the case that $m > N \geq k$, we now show how $M_{k}$ can be constructed in $O(kmn)$ flops.


First, we note that $M_0$ can be constructed explicitly using \cref{M0} in $\mO(m_0n)$ flops. Suppose that $(a_0, a_1, \ldots)$ is the infinite coefficient vector of the Chebyshev approximant to a univariate function $a(x)$. It has been shown that $\mM_1[a(x)]$ is Toeplitz-plus-Hankel \cite{olv1,qin1}
\begin{align}
\mM_1[a(x)] = \frac{1}{2}\left[
\begin{pmatrix}
2a_0 & a_1 & a_2 & a_3 & \cdots \\
a_1 & 2a_0 & a_1 & a_2 & \ddots \\
a_2 & a_1 & 2a_0 & a_1 & \ddots \\
a_3 & a_2 & a_1 & 2a_0 & \ddots \\
\vdots & \ddots & \ddots & \ddots & \ddots \\
\end{pmatrix} -
\begin{pmatrix}
a_2 & a_3 & a_4 & a_5 & \cdots \\
a_3 & a_4 & a_5 & a_6 & \udots \\
a_4 & a_5 & a_6 & a_7 & \udots \\
a_5 & a_6 & a_7 & a_8 & \udots \\
\vdots & \udots & \udots & \udots & \udots \\
\end{pmatrix}\right]. \label{M1structure}
\end{align}
Therefore, $M_1$ can be explicitly constructed in $\mO(m_1n)$ flops. For $k \geq 2$, the construction of $M_{k}$ is done via a detour to $M_1$, as implied by the following lemma.
\begin{lemma}
For a univariate function $a(x)$ and $k \geq 2$, the multiplication operator $\mM_{k}[a(x)]$ can be represented as
\begin{align}
\mM_{k}[a(x)] = \mS_{k - 1}\mS_{k - 2}\ldots \mS_{1} \mM_{1}[a(x)]\mS_{1}^{-1} \mS_{2}^{-1} \ldots \mS_{k - 1}^{-1}.\label{Mlambda}
\end{align}
In addition, if $a(x)$ is a finite Chebyshev series of degree $m_{k}$, $\mM_1[a(x)]$ and $\mM_{k}[a(x)]$ are both banded with lower and upper bandwidths $(m_{k}, m_{k})$.
\end{lemma}

\begin{proof}
Suppose that $u(x) = \sum_{j = 0}^{\infty} u_j^{k-1} C_{j}^{(k-1)}(x)=\sum_{j = 0}^{\infty} u_j^{k} C_{j}^{(k)}(x)$, $\bsu^{k-1}=(u_0^{k-1},$ $u_1^{k-1}, \ldots)^{\top}$, and $\bsu^{k}=(u_0^{k}, u_1^{k}, \ldots)^{\top}$. By the definition of $\mM_{k}$ we have
\begin{align*}
  \left(C_0^{(k)} \,\middle|\, C_1^{(k)} \,\middle|\, \cdots \right) \mM_{k}[a(x)] \bsu^{k} = 
  \left(C_0^{(k - 1)} \,\middle|\, C_1^{(k - 1)} \,\middle|\, \cdots \right) \mM_{k - 1}[a(x)] \bsu^{k - 1},
\end{align*}  
as both sides represent $a(x)u(x)$. Since $\mS_{k - 1}$ transforms the coefficients in $C^{(k - 1)}$ to those in $C^{(k)}$,
\begin{align*}
  \left(C_0^{(k - 1)} \,\middle|\, C_1^{(k - 1)} \,\middle|\, \cdots \right) \mM_{k - 1}[a(x)] \bsu^{k - 1} = 
  \left(C_0^{(k)} \,\middle|\, C_1^{(k)} \,\middle|\, \cdots \right) \mS_{k - 1} \mM_{k-1}[a(x)] \mS_{k - 1}^{-1} \bsu^{k}.
\end{align*}  
\cref{Mlambda} follows from induction. 

The assertion on the bandwidths follows from the three-term recurrence relation \cref{M}. See \cite[\S 6.3.1]{tow1} for details.
\end{proof}

In light of \cref{Mlambda}, $M_{k}[a^{k}]$ can be approximated as $S_{k-1}\dots S_1 M_1[a^{k}]S^{-1}_1 \dots S^{-1}_{k-1}$. Suppose that we have constructed $M_1[a^{k}]$ using \cref{M1structure}. Let $Y = S_1 M_1[a^{k}]$ and $X = Y S_1^{-1}$. Since the bandwidths of $X$ are also $(m_{k}, m_{k})$, we first compute $Y$ by only calculating the entries in the $(m_{k}, m_{k})$ band, followed by solving
\begin{align*}
S_1^T X^T = Y^T
\end{align*}
for only the entries of $X$ in the $(m_{k}, m_{k})$ band. In fact, $Y$ is banded with bandwidths $(m_{k}, m_{k}+2)$, but the entries in the $(m_{k}+1)$th and $(m_{k}+2)$th superdiagonals are not involved in the calculation of the entries of $X$ in the $(m_{k}, m_{k})$ band. Hence, the similarity transform of $M_1[a^{k}]$ with respect to $S_1$ costs $\mO(m_{k}n)$ flops. We proceed in the same manner for the similarity transforms with $S_2, S_3, \ldots, S_{k - 1}$, entailing a cost of $\mO(k m_{k} n)$ for constructing a specific $M_{k}[a^{k}]$. It then follows that the total cost for forming all $M_{k}[a^{k}]$ is $\mO(N^2 m n)$. 

\begin{figure}[t!]
  \centering
  \subfloat[$k = 10$ and $n = 10^4$]{\label{fig:speed_mul_c10}\includegraphics[scale=0.45]{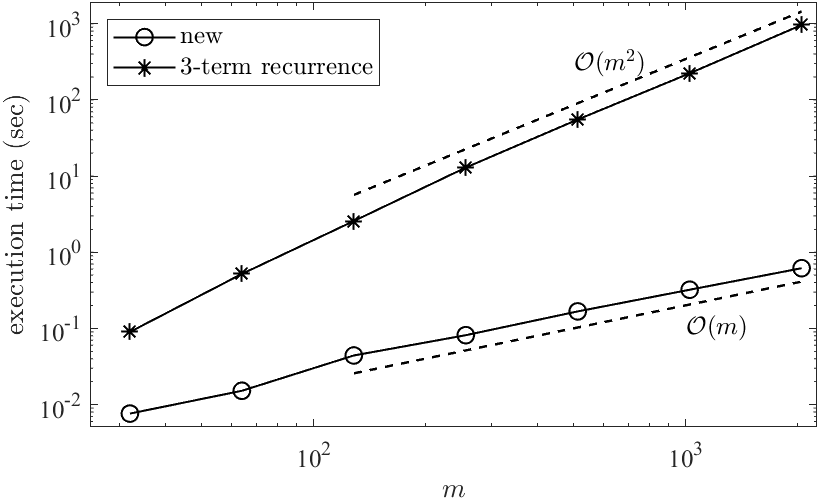}}
  \hfill
  \subfloat[$k = 2$ and $m = 20$]{\label{fig:speed_mul_c2}\includegraphics[scale=0.45]{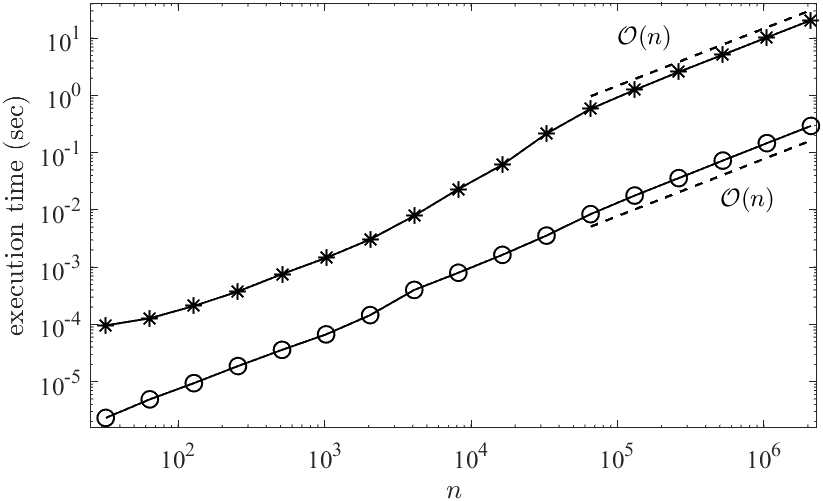}}
  \caption{Execution times for constructing $M_{k}[a(x)]$ using \cref{Mlambda} and the method described in \cite[\S 6.3.1]{tow1}.}\label{fig:speed_mul}
\end{figure}

In \cref{fig:speed_mul}, the proposed approach to constructing $M_{k}$ is compared against \cref{M}. We first have $k$ and $n$ fixed ($k=10$ and $n=10^4$) but increment $m$ to examine the dependence on $m$. As $m$ varies from $2^5$ to $2^{11}$, the proposed method is approximately $11\times$ to $1563\times$ faster than the recurrence method. In \cref{fig:speed_mul_c2}, we let $k=2$ and $m=20$ and vary $n$ from $2^5$ and $2^{21}$. Now the two methods have the same asymptotics, but the proposed method is still much faster---the speedup ranges between $20\times$ to $70\times$.

Instead of assembling $L$ plainly as suggested by \cref{L}, we construct $L$ in a nested fashion by imitating the Horner's method
\begin{align}
L = M_N[a^N]D_N + S_{N-1} \biggl(\cdots + S_{2} \Bigl( M_2[a^2]D_2 +S_1 \left(M_1[a^1] D_1 + S_0M_0[a^0] \right) \Bigr) \biggr). \label{Lnested}
\end{align}
Although it still incurs $\mO(N^2mn)$ flops, \cref{Lnested} saves half of the cost from plain calculation. The overall complexity for obtaining $L$ is $\mO(N^2mn)$. Since it is usually that $m>N$ or $m\gg N$, this is in contrast to $\mO(Nm^2n)$, the complexity of the recurrence method or the explicit formula.

Finally, forming the product in \cref{coe_banded_pg} incurs another cost of $\mO(Nmn)$ flops. Thus, the overall cost of constructing $A$ is $\mO(N^2mn)$.

\section{Numerical examples} \label{sec:exm}
In this section, we demonstrate the efficiency and the accuracy of the new banded PG spectral method by two examples. We compare the new method against MPG method \cite{mor1} and the US method \cite{olv1}. For MPG method, the coefficient matrix on the left-hand side can be constructed using two different approaches. If the variable coefficients of the ODE all can be written as power series, a recursion (R) approach can be taken at a cost of $\mO(n)$ to construct the system directly. Otherwise, one has to resort to numerical integration (NI) for evaluating the entries of the coefficient matrix. The construction of the linear system for the US method follows \cite{olv1}, except that the multiplication matrices $M_{k}$ are constructed using the three-term recurrence method \cref{M}.

Suppose that $u(x)$ and $\hat{u}(x)$ are the exact and a computed solution respectively. We measure the error using the absolute $2$-norm
$\left(\int_{-1}^{1}(u(x) - \hat{u}(x))^2 \mathrm{d} x \right)^{1/2}$.

All the experiments are conducted in \textsc{Julia} v1.10.2 on a laptop with a 4-core 2.8 GHz Intel i7-1165G7 CPU. Execution times are measured using \texttt{Benchmark.jl}.

\subsection{An ODE with simple variable coefficients}
\begin{figure}[t!]
\centering
\subfloat[construction time]{\label{fig:time_taylor}\includegraphics[scale=0.45]{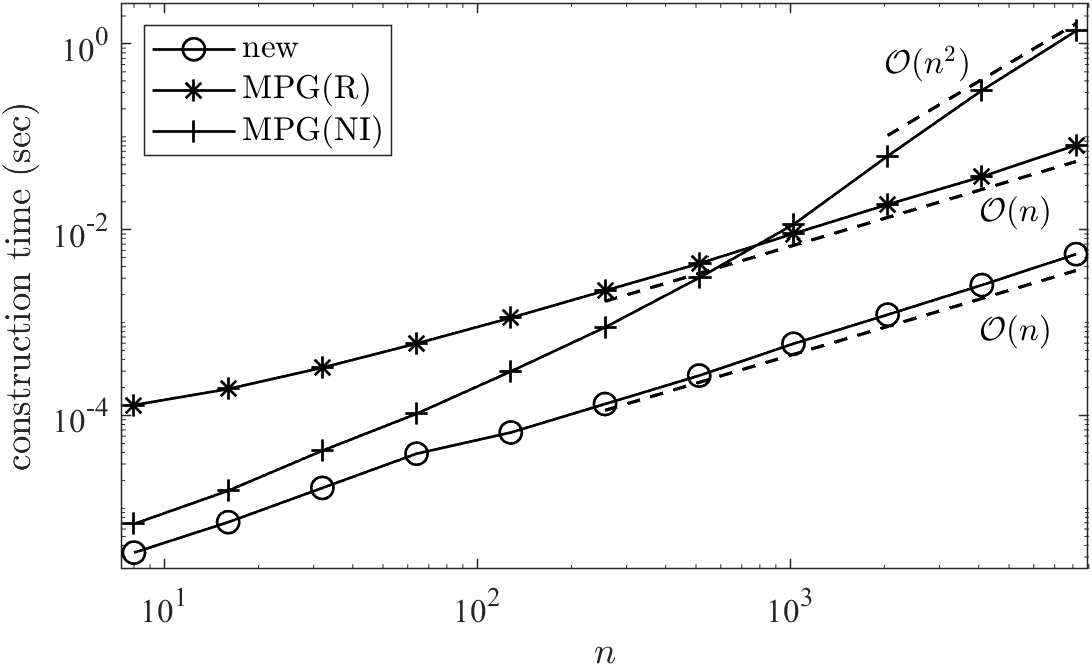}}
\hfill
\subfloat[absolute error]{\label{fig:err_taylor}\includegraphics[scale=0.45]{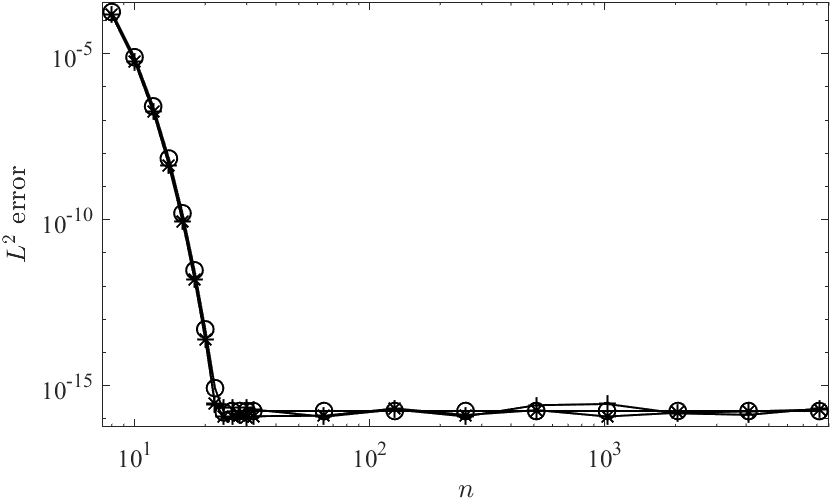}}
\caption{Solving \cref{taylor} with the proposed and MPG methods.}\label{fig:taylor}
\end{figure}
Our first example, adapted from \cite{she5}, is
\begin{align}
u''' - \cos(x)u'' + 10e^{x}u = f(x), ~ u(\pm 1) = 1, ~ u'(1) = 0  \label{taylor}
\end{align}
with $f(x)$ chosen so that $u(x) = e^{(x^2 - 1) / 2}$. We use this example to demonstrate the advantage of the proposed method over MPG method on the construction of the linear system. Since both MPG and the proposed method lead to strictly banded matrices with similar bandwidths, we omit the comparison of the solution times as they are almost the same. The variable coefficients of this ODE have known Taylor expansions. Thus, MPG method can construct the linear system via recursion (see \cite[Table 1]{mor1}) and only incurs a cost proportional to $n$. Since \cref{taylor} is of an odd order, the test functions are chosen to satisfy the dual boundary conditions. 

\cref{fig:time_taylor} displays the execution time on constructing the linear system taken by MPG(R), MPG(NI), and the proposed method for various $n$. For MPG(NI) method, each of the $\mO(n)$ nonzero entries in the band incurs a cost of $\mO(n)$ flops if it is evaluated by numerical integration. This matches the curve for MPG(NI) method, which grows quadratically. The curves for MPG(R) and the proposed method both exhibit linear complexity in construction time. However, as anticipated, the proposed method is $14\times$ faster. We also note that it is until about $n = 10^3$ MPG(R) method is faster than MPG(NI), implying a large hidden constant in the big-Oh notation for MPG(R) method. Since the recursion technique is applicable only when the variable coefficients are sufficiently simple, the cost of numerical integration limits the practical performance of MPG in term of speed, especially when large discretization size is required. \cref{fig:err_taylor} shows that the convergences of MPG(R), MPG(NI), and the proposed method behave similarly.

\begin{figure}[t!]
  \centering
  \subfloat[execution time]{\label{fig:time_airy}\includegraphics[scale=0.45]{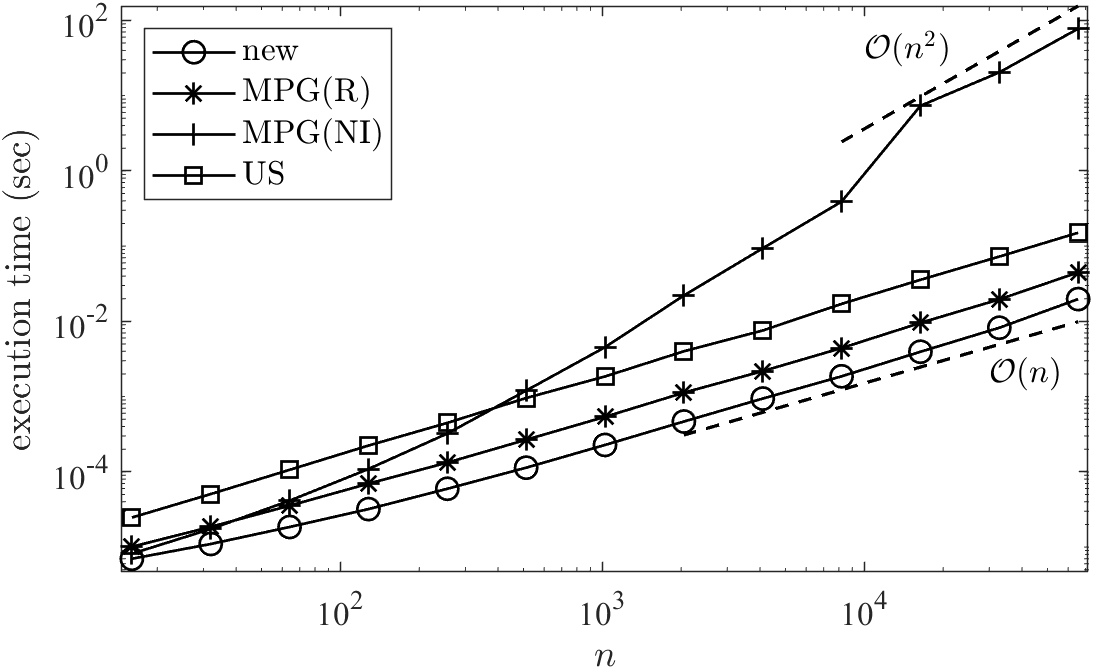}}
  \hfill
  \subfloat[absolute error]{\label{fig:err_airy}\includegraphics[scale=0.45]{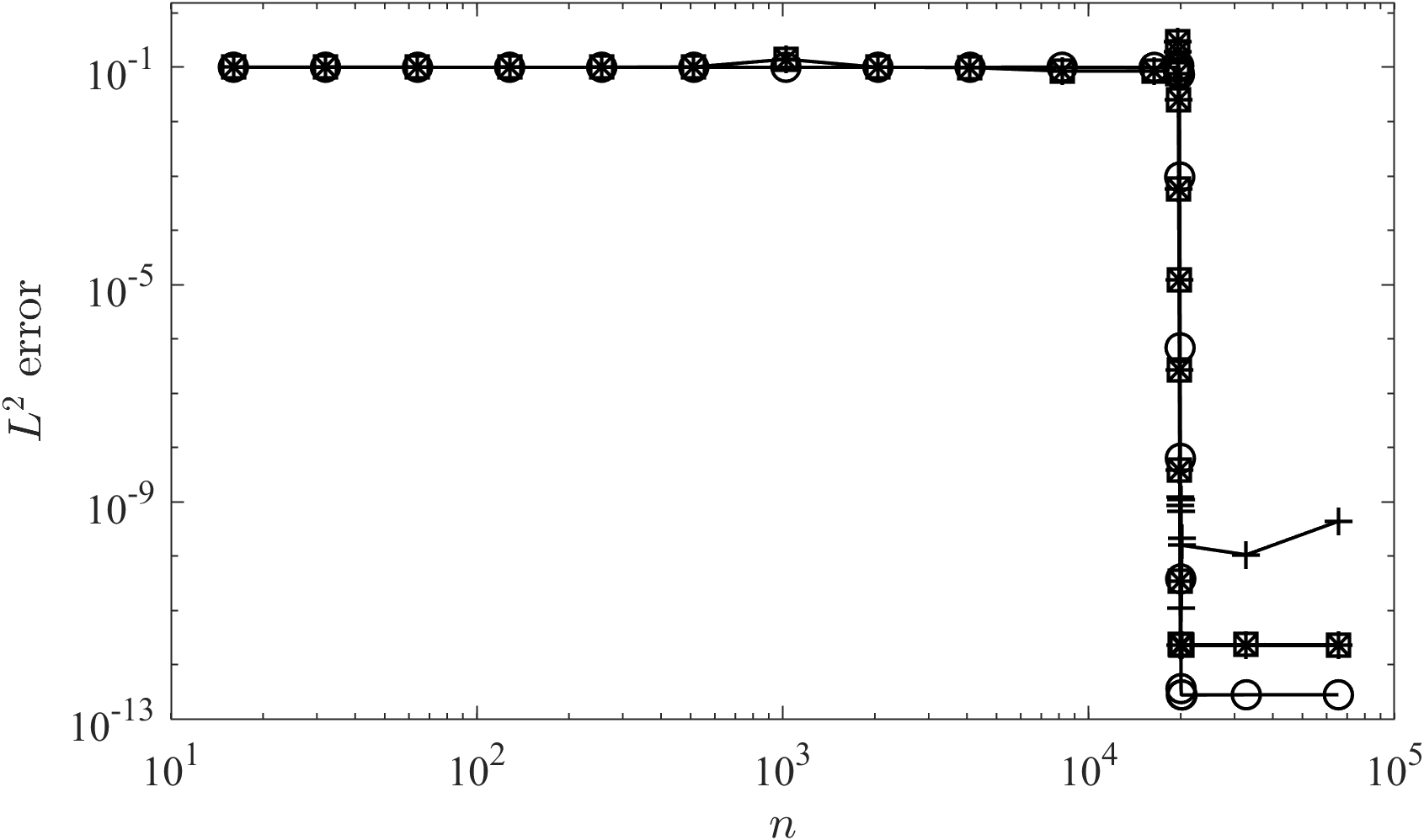}}
  \caption{Solving \cref{airy} for $\epsilon = 10^{-9}$ using the proposed, MPG(R), MPG(NI), and the US methods.}\label{fig:airy}
\end{figure}

\subsection{Airy Equation}
Our second example is the Airy equation in the canonical domain
\begin{align}
\epsilon u'' - x u(x) = 0,~~~
u(-1) = \operatorname{Ai}\left(-\sqrt[3]{\epsilon^{-1}} \right),~~~
u(1) = \operatorname{Ai}\left(\sqrt[3]{\epsilon^{-1}} \right), \label{airy}
\end{align}
where $\operatorname{Ai}(x)$ is the Airy function of the first kind. The exact solution is $\operatorname{Ai}(\sqrt[3]{\epsilon^{-1}}x)$. We solve \cref{airy} for $\epsilon = 10^{-9}$. For such a small $\epsilon$, \cref{airy} features a lengthy solution, serving as an ideal test problem for solution speed. Apparently, MPG(R) works readily for this problem. Besides MPG and the proposed method, we also include the US method for comparison. The almost-banded system arising from the US method is solved by the \texttt{qr} function in \texttt{SemiseparableMatrices.jl}. For the proposed method, the test space is chosen to be the same as the trial.

\cref{fig:time_airy} shows the total execution time including the construction and the solution. Although both the construction and the solution have linear complexities for MPG(R), the US, and the proposed method, the new method is at least $2\times$ and $8\times$ as fast as MPG(R) and the US method respectively. This is because the new method is faster in construction and a banded system, when solved by the \textsc{Lapack} routines such as \texttt{gbtrf} and \texttt{gbtrs}, has a clear edge of speed over an almost banded system that can only be solved by user-supplied code that is unlikely to be optimized in terms of memory caching and allocation.

As shown in \cref{fig:err_airy}, spectral convergence takes place for all the methods at somewhere between $n=19,500$ and $n=20,100$. The final accuracy of MPG(NI) method is $O(10^{-10})$, two orders lower than those of MPG(R) and the US methods. Among the four methods, the proposed method is most accurate with an error of $O(10^{-13})$. In fact, the entrywise difference of the coefficient matrices between MPG(NI) and MPG(R) methods is less than $O(10^{-17})$ and the much magnified discrepancy in the accuracy is the consequence of the poor conditioning of the problem. This example, along with our extensive experiments, implies that numerical integration is less favorable for singularly perturbed problems due to the tampered accuracy. Thus, the fact that the MPG method often relies on numerical integration for general variable-coefficient ODEs may somewhat limit its applicability as a general ODE solver. 

\section{An overarching PG method} \label{sec:jacobi}
In this section, we generalize the new PG method to Jacobi polynomials $J_{n}^{(\alpha, \beta)}(x)$ and, much more importantly, show that MPG method is a particular case of this generalization. Since existing banded Galerkin spectral methods are specific instances of MPG method, as shown in \cite{mor1}, our Jacobi-based generalization can be viewed as the overarching method for all existing banded Galerkin spectral methods. For notational convenience, we denote by $\omega^{(\alpha, \beta)}(x)$ the Jacobi weight function $(1-x)^{\alpha}(1+x)^{\beta}$. 

\subsection{Jacobi-based banded PG method}
We start off by first setting out the Jacobi-based operators that are analogous to those used by the US method. The properties of Jacobi polynomials that we use below can be found in standard texts, e.g., \cite[\S 8.2]{luk1}. 

For $k \geq 1$, the Jacobi-based $k$th differential operator
\begin{align*} 
  \begin{aligned}
    & ~~\,\,\, \myunderbrace{\text{\footnotesize $k$ times}} \\[-9pt]
    \mD_{k}^{J} = \frac{1}{2^k} & \begin{pmatrix}
         \  0 ~~  \cdots ~~ 0  &\frac{\Gamma(\alpha + \beta + 2k + 1)}{\Gamma(\alpha + \beta + k + 1)} & & & \\
        & &\frac{\Gamma(\alpha + \beta + 2k + 2)}{\Gamma(\alpha + \beta + k + 2)}& & \\
        & & &\frac{\Gamma(\alpha + \beta + 2k + 3)}{\Gamma(\alpha + \beta + k + 3)}\\
        & & & & \ddots\\
    \end{pmatrix}
  \end{aligned},
\end{align*}
which satisfies
\begin{align*}
  \frac{\mathrm{d}^k}{\mathrm{d}x^k}\left(J_{0}^{(\alpha, \beta)} \,\middle|\, J_{1}^{(\alpha, \beta)} \,\middle|\, \cdots \right) 
  = \left(J_{0}^{(\alpha+k, \beta+k)} \,\middle|\, J_{1}^{(\alpha+k, \beta+k)} \,\middle|\, \cdots \right)\mD_{k}^{J}.
\end{align*} 

For a variable coefficient, we need the operator $\mM_{k}^{J}[h(x)]$ to represent the pre-multiplication by a given Jacobi series, say, e.g., $h(x)$. That is,
\begin{align*}
  h(x) \left(J_{0}^{(\alpha+k, \beta+k)} \,\middle|\, J_{1}^{(\alpha+k, \beta+k)} \,\middle|\, \cdots \right) = \left(J_{0}^{(\alpha+k, \beta+k)} \,\middle|\, J_{1}^{(\alpha+k, \beta+k)} \,\middle|\, \cdots \right) \mM_{k}^{J}[h(x)].
\end{align*} 
Since no closed-form expression is known, $\mM_{k}^{J}[h(x)]$ can only be constructed via the three-term recurrence relation
\begin{align*}
\mM_{k}^{J}[J_{j + 1}^{(\alpha+k, \beta+k)}] &= \left(a_{j}^{(\alpha+k, \beta+k)}\mM_{k}^{J}[x] - b_{j}^{(\alpha+k, \beta+k)}\mI\right)\\
&\times \mM_{k}^{J}[J_{j}^{(\alpha+k, \beta+k)}] - c_{j}^{(\alpha+k, \beta+k)}\mM_{k}^{J}[J_{j - 1}^{(\alpha+k, \beta+k)}], \quad j \geq 1, \label{MJ}
\end{align*}
where the recurrence coefficients
\begin{align*}
a_{j}^{(\alpha+k, \beta+k)} &= \frac{(2j + \alpha + \beta + 2k + 1)(2j + \alpha + \beta + 2k + 2)}{2(j + 1)(j + \alpha + \beta + 2k+ 1)}, \\ 
b_{j}^{(\alpha+k, \beta+k)} &= \frac{((\beta+k)^2 - (\alpha+k)^2)(2j + \alpha + \beta + 2k + 1)}{2(j+1)(j + \alpha + \beta + 2k + 1)(2j + \alpha + \beta + 2k)},\\ 
c_{j}^{(\alpha+k, \beta+k)} &= \frac{(j + \alpha + k)(j + \beta + k)(2j + \alpha + \beta + 2k + 2)}{(j + 1)(j + \alpha + \beta + 2k + 1)(2j + \alpha + \beta + 2k)}.
\end{align*}
The recursion is started off with $\mM_{k}^{J}[J_{0}^{(\alpha+k, \beta+k)}]=\mI$, $\mM_{k}^{J}[J^{(\alpha+k, \beta+k)}_{1}]=a_{0}^{(\alpha+k, \beta+k)}\mM_{k}^{J}[x] - b_{0}^{(\alpha+k, \beta+k)}\mI$, and
\begin{align*}
  \mM_{k}^{J}[x] =
  \begin{pmatrix}
    \tilde{b}_{0}^{(\alpha+k, \beta+k)} & \tilde{c}_{1}^{(\alpha+k, \beta+k)} &  &  &  \\
    \tilde{a}_{0}^{(\alpha+k, \beta+k)} & \tilde{b}_{1}^{(\alpha+k, \beta+k)} & \tilde{c}_{2}^{(\alpha+k, \beta+k)} &  &  \\
     & \tilde{a}_{1}^{(\alpha+k, \beta+k)} & \tilde{b}_{2}^{(\alpha+k, \beta+k)} & \tilde{c}_{3}^{(\alpha+k, \beta+k)} &  \\
     &  & \tilde{a}_{2}^{(\alpha+k, \beta+k)} & \tilde{b}_{3}^{(\alpha+k, \beta+k)} & \tilde{c}_{4}^{(\alpha+k, \beta+k)} \\
     &  &  \ddots & \ddots  & \ddots
  \end{pmatrix} ~\text{for }k \geq 0, 
\end{align*}  
where
\begin{align*}
  \tilde{a}_{j}^{(\alpha+k, \beta+k)} &= \frac{2(j + 1)(j + \alpha + \beta + 2k + 1)}{(2j + \alpha + \beta + 2k + 1)(2j + \alpha + \beta + 2k + 2)},\\
  \tilde{b}_{j}^{(\alpha+k, \beta+k)} &= \frac{(\beta + k)^2 - (\alpha + k)^2}{(2j + \alpha + \beta + 2k)(2j + \alpha + \beta + 2k + 2)},
  \\ 
  \tilde{c}_{j}^{(\alpha+k, \beta+k)} &= \frac{2(j + \alpha + k)(j + \beta + k)}{(2j + \alpha + \beta + 2k)(2j + \alpha + \beta + 2k + 1)}.
\end{align*}
For $h(x) = \sum_{j = 0}^{d}h_jJ_{j}^{(\alpha+k, \beta+k)}(x)$, it is straightforward to show that $\mM_{k}^{J}[h(x)] = \sum_{j = 0}^{d}h_j\mM_{k}^{J}[J_{j}^{(\alpha+k, \beta+k)}]$ is a banded matrix with bandwidths $(d, d)$. 

Like in the ultraspherical case, we need the conversion operator 
\begin{align*}
  \mS_{k}^{J} = \begin{pmatrix}
    \ha_{0}^{(\alpha+k, \beta+k)}&\hb_{1}^{(\alpha+k, \beta+k)}& \hc_{2}^{(\alpha+k, \beta+k)}& & \\ 
   &\ha_{1}^{(\alpha+k, \beta+k)}&\hb_{2}^{(\alpha+k, \beta+k)}&\hc_{3}^{(\alpha+k, \beta+k)}& \\
   & &\ha_{2}^{(\alpha+k, \beta+k)}&\hb_{3}^{(\alpha+k, \beta+k)}&\hc_{4}^{(\alpha+k, \beta+k)}\\
   & & \ddots&\ddots&\ddots\\
  \end{pmatrix} ~\text{for } k \geq 0,
\end{align*}
such that
\begin{align*}
  \left(J_{0}^{(\alpha+k, \beta+k)} \,\middle|\, J_{1}^{(\alpha+k, \beta+k)} \,\middle|\, \cdots \right) 
  = \left(J_{0}^{(\alpha+k+1, \beta+k+1)} \,\middle|\, J_{1}^{(\alpha+k+1, \beta+k+1)} \,\middle|\, \cdots \right)\mS_{k}^{J}.
\end{align*}  
The nonzero entries in $\mS_{k}^{J}$ are
\begin{align*}
  \ha_{j}^{(\alpha+k, \beta+k)} &= \frac{(j + \alpha + \beta + 2k + 1)(j + \alpha + \beta + 2k + 2)}{(2j + \alpha + \beta + 2k + 1)(2j + \alpha + \beta + 2k + 2)}, \\ 
  \hb_{j}^{(\alpha+k, \beta+k)} &= \frac{(\alpha - \beta)(j + \alpha + \beta + 2k + 1)}{(2j + \alpha + \beta + 2k)(2j + \alpha + \beta + 2k +2)}, \\ 
  \hc_{j}^{(\alpha+k, \beta+k)} &= -\frac{(j+\alpha + k)(j + \beta + k)}{(2j + \alpha + \beta + 2k)(2j + \alpha + \beta + 2k + 1)},
  \end{align*}
for $j \geq 0$. The three operators that we have just spelled out can be used to construct the Jacobi-based US method, although this is not the goal here and nothing is gained from involving Jacobi polynomials for the US method.

Since the trial and test functions $\phi^{J}_{k}(x)$ and $\psi^{J}_k(x)$ are recombined Jacobi polynomials, they are written as
\begin{align*}
  \left(\phi^{J}_0 \,\middle|\, \phi^{J}_1 \,\middle|\, \cdots \,\middle|\, \phi^{J}_{n-1} \right) 
  &= \left(J^{(\alpha, \beta)}_0 \,\middle|\, J^{(\alpha, \beta)}_1 \,\middle|\, \cdots \,\middle|\, J^{(\alpha, \beta)}_{n+N-1} \right) R^{J}, \\
  \left(\psi^{J}_0 \,\middle|\, \psi^{J}_1 \,\middle|\, \cdots \,\middle|\, \psi^{J}_{n-1} \right)
  &= \left(J^{(\alpha + N, \beta+N)}_0 \,\middle|\, J^{(\alpha + N, \beta+N)}_1 \,\middle|\, \cdots \,\middle|\, J^{(\alpha + N, \beta+N)}_{n+N-1} \right) Q^{J},
\end{align*}  
where $Q^{J}$ and $R^{J}$ are the stencil matrices. 

Now we have all the ingredients for constructing the Jacobi-based banded PG spectral method as stated in the following lemma. We omit the proof, for it is analogous to \cref{lem:banded_cmtx}. Particularly, setting $\alpha = \beta = -1/2$ reduces \cref{lem:bjpg} to \cref{lem:banded_cmtx}, up to a scaling factor.

\begin{lemma} \label{lem:bjpg}
For the PG spectral method with trial functions $\phi^{J}_{k}(x)$, test functions $\psi^{J}_{k}(x)$, and the positive weight function $\omega^{(\alpha + N, \beta + N)}(x)$, the coefficient matrix
\begin{align}
A^{J} = (Q^{J})^{T} \Omega^J_{n + N - 1} L^J_{n + N - 1} R^{J}, \label{Ajacobi}
\end{align}
where
\begin{align*}
  \Omega^{J}_{n + N - 1} = \diag \left(d_0^J, \ldots, d_{n+N-1}^{J}\right), \; 
  d_{j}^{J} = \frac{2^{\alpha+\beta+2N+1}\Gamma(j+\alpha+N+1)\Gamma(j+\beta+N+1)}{\left(2j+\alpha+\beta+2N+1\right)j!\Gamma(j+\alpha+\beta+2N+1)}
\end{align*}
and 
\begin{align*}
L^{J}_{n + N - 1} = \mP_{n+N-1}\left(\mM_N^{J}[a^N] \mD^{J}_N + \sum_{k = 0}^{N-1}{\mS}^{J}_{N-1}\ldots {\mS}^{J}_{k}\mM^{J}_{k}[a^{k}]\mD^{J}_{k}\right)\mP_{n+N-1}^{\top},
\end{align*}
where $\mD_0^J = \mI$.
\end{lemma}

\subsection{MPG method as a specific instance}
We now turn to the main goal of this section---show that MPG method is a specific instance of the new Jacobi-based PG spectral method. To see this, we note that MPG method uses exactly the same combinations of Jacobi polynomial $J_{k}^{(\alpha, \beta)}$ as trial functions, whereas it takes $\left\{ \left(h^{(N)}_{N+k}\right)^{-1}\left(1-x^2 \right)^N  \partial^{N} J_{k + N}^{(\alpha, \beta)} \right\}_{k = 0}^{n-1}$ as the test functions. Here, $h^{(N)}_{N+k}$ is a normalization factor for the $k$th test function (see \cite[Equation (2.16)]{mor1}). With the weight function $\omega^{(\alpha, \beta)}$, the MPG coefficient matrix
\begin{align*}
  A^{MPG} = H^{-1}
  \left(
  \begin{pmatrix}
            (1-x^2)^N\partial^{N}J_{N}^{(\alpha, \beta)}\\
            (1-x^2)^N\partial^{N}J_{N + 1}^{(\alpha, \beta)}\\
            \vdots\\
            (1-x^2)^N\partial^{N}J_{n + N - 1}^{(\alpha, \beta)}
           \end{pmatrix}
           \mL \left(J_{0}^{(\alpha, \beta)} \,\middle|\, J_{1}^{(\alpha, \beta)} \,\middle|\, \cdots \,\middle|\, J_{n + N - 1}^{(\alpha, \beta)} \right)
    \right)_{\omega^{(\alpha, \beta)}} \hspace{-5mm} R^{J}, \\ 
\end{align*}  
where $H = \operatorname{diag}\left(h^{(N)}_N, h^{(N)}_{N+1}, \dots, h^{(N)}_{n+ N -1}\right)$. Further, moving the factor $(1-x^2)^N$ to the weight function gives
\begin{align}
  A^{MPG} &= H^{-1}\left(
      \begin{pmatrix}
            \partial^{N}J_{N}^{(\alpha, \beta)}\\
            \partial^{N}J_{N + 1}^{(\alpha, \beta)}\\
            \vdots\\
            \partial^{N}J_{n + N - 1}^{(\alpha, \beta)}
           \end{pmatrix}
           \mL \left(J_{0}^{(\alpha, \beta)} \,\middle|\, J_{1}^{(\alpha, \beta)} \,\middle|\, \cdots \,\middle|\, J_{n + N - 1}^{(\alpha, \beta)} \right) 
    \right)_{\omega^{(\alpha + N, \beta + N)}} \hspace{-13mm} R^{J} \nonumber \\ 
    &= H^{-1} (D_{N}^{J})^{T} \left(
      \begin{pmatrix}
            J_{0}^{(\alpha + N, \beta + N)}\\
            J_{1}^{(\alpha + N, \beta + N)}\\
            \vdots\\
            J_{n + N - 1}^{(\alpha + N, \beta + N)}
           \end{pmatrix}
           \left(J_{0}^{(\alpha + N, \beta + N)} \,\middle|\, J_{1}^{(\alpha + N, \beta + N)} \,\middle|\, \cdots \,\middle|\, J_{n + N - 1}^{(\alpha + N, \beta + N)} \right) L_{n + N - 1}^{J} 
    \right)_{\omega^{(\alpha + N, \beta + N)}} \hspace{-13mm} R^{J} \nonumber \\
    &= H^{-1} (D_{N}^{J})^{T}\Omega_{n + N - 1}^{J}L_{n + N - 1}^{J} R^{J}, \label{Ampg}
\end{align}
where $D^J_N = \mP_{n+N-1}\mD_{N}^J \tilde{\mP}_{n}^{\top}$ and $\tilde{\mP}_{n} = \left(0_{n \times N}, \mP_{n}\right)$ is a shifted projection operator. If we take $(Q^{J})^{T} = H^{-1}(D_{N}^{J})^{T} $, \cref{Ajacobi} and \cref{Ampg} become identical. 

The fact that the new method encompasses the MPG method and other historically prominent banded spectral Galerkin methods as specific instances positions it as a lens through which existing methods can be examined and studied from a broader perspective. Furthermore, it suggests that designing new sparse Galerkin spectral methods following traditional patterns may lead to approaches that fit within the proposed framework.


\section{Accelerating the US method}\label{sec:us}
\begin{figure}[t!]
  \centering
  \subfloat[execution time for construction]{\label{fig:tenth_con}\includegraphics[scale=0.45]{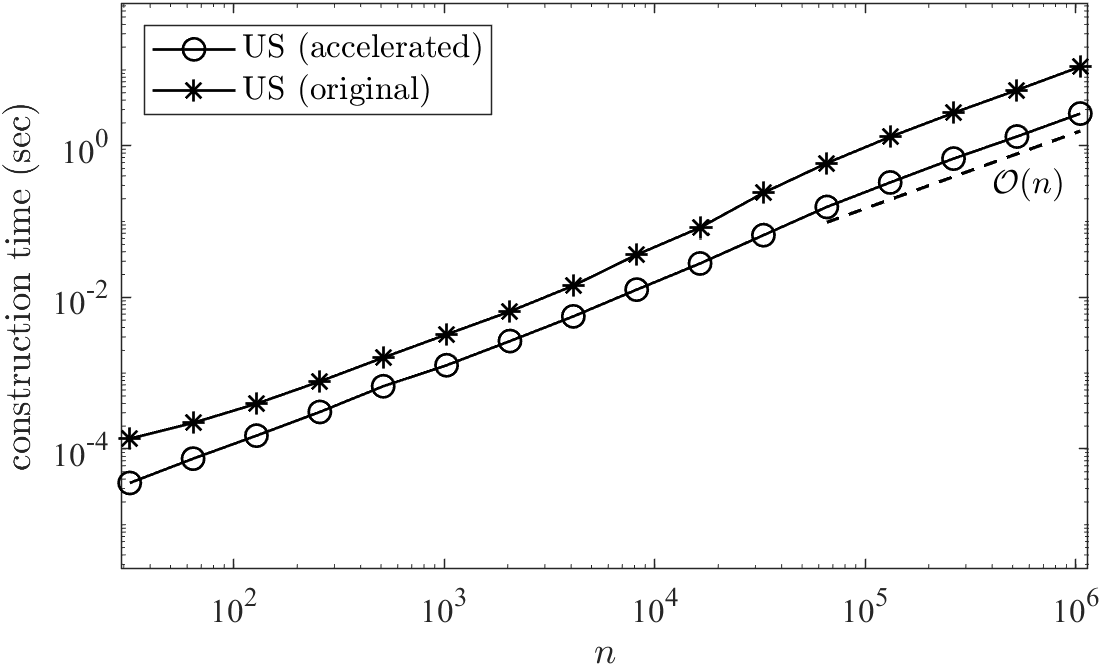}}
  \hfill
  \subfloat[execution time for solution]{\label{fig:tenth_sol}\includegraphics[scale=0.45]{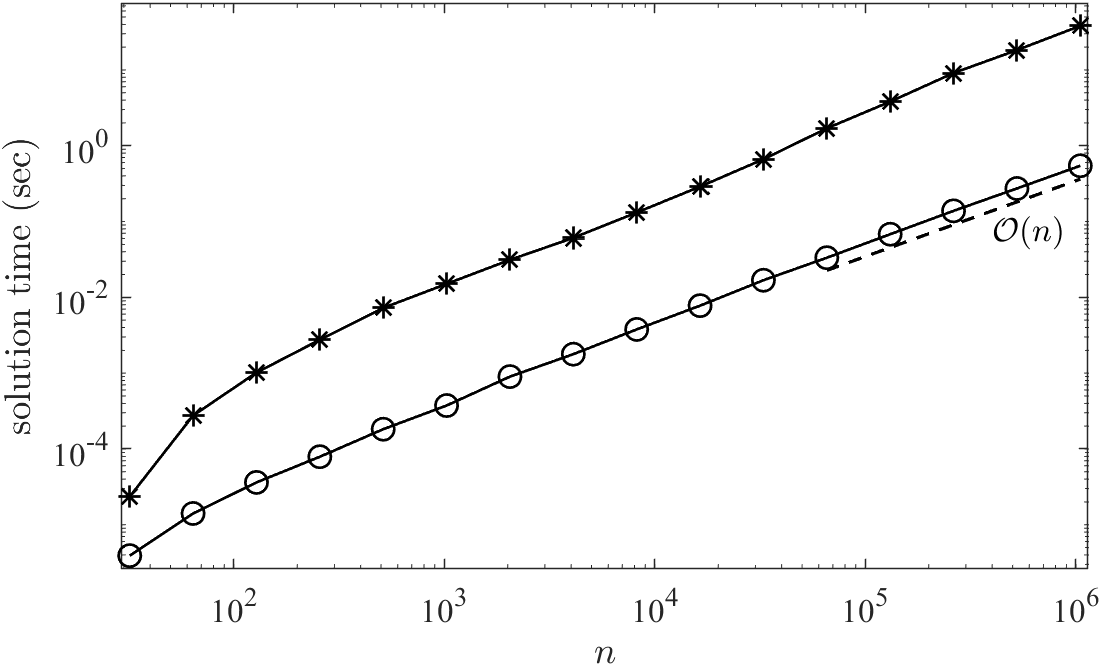}}
  \caption{Solving the $10$th-order ODE \cref{tenth} using the original and the accelerated US methods.}\label{fig:tenth}
\end{figure}
Nothing holds us from accelerating the US method by the techniques introduced in \cref{sec:bpg,sec:fs}. When doing so, we actually migrate from a tau method to a PG approach. Particularly, in case of homogeneous linear constraints the accelerated US method becomes an instance of the framework proposed by \cref{lem:banded_cmtx}, if $Q$ is set to a rectangular truncation of the identity matrix. Consider the $10$th-order ODE
\begin{align}
  \label{tenth}
  \begin{aligned}
    &u^{(10)} + \cosh(x) u^{(8)} + x^2 u^{(6)} + x^4 u^{(4)} + \cos(x) u^{(2)} + x^2 u = 0,\\
    &u'(-1)=u'(1)=1 \text{ and } u^{(k)}(\pm 1) = 0 \text{ for } k =0,2,3,4,
  \end{aligned}
\end{align}
which is taken from \cite{olv1}. We accelerate the original US method by basis recombination, fast construction of the multiplication operators, and the nested assembly and show in \cref{fig:tenth} the comparison with the non-accelerated version. The speedups in construction and solution for $n$ up to $10^6$ are at least $2\times$ and $40\times$ respectively.

\section{Discussion}\label{sec:closing}
Basis recombination, as a strategy of enforcing boundary conditions or other side constraints, is often compared to boundary bordering, which is used in tau method and the US method. It is argued that the downsides of basis recombination are as follows.
\begin{enumerate}
\item The solution is not computed in the convenient and orthogonal basis.
\item Without orthogonal polynomials one cannot apply recurrence relations to construct multiplication matrices as \cref{M}.
\item The structure of the linear systems may depend on the boundary conditions or side constraints, therefore prohibiting the use of a fast, general solver.
\item The solution may be expressed in an unstable basis for problems with very high-order boundary conditions.
\item There is no unique way of combining an orthogonal basis.
\end{enumerate}

The current investigation however leads to somewhat different observations. First, although the solution is obtained in the recombined basis, transforming it back to the coefficients in the orthogonal basis is straightforward---simply premultiply the solution vector by $R$. Second, constructing the multiplication matrices directly using the recurrence relation is shown to be slow. For Chebyshev and ultraspherical polynomials, we recommend the fast construction introduced in \cref{sec:fs}, particularly when the variable coefficients of the ODE can only be approximated by Chebyshev series of large degrees. Third, as we have shown, the resulting systems are always banded; they only differ in their bandwidths. Nonetheless, this has little effect on the performance of banded solvers. Fourth, even for high-order boundary conditions, the recombined basis is numerically stable due to the safeguard value $\gamma_{\min}$. Fifth, we have shown in \cref{sec:tt} that if we aim at minimal bandwidths in the resulting systems the way of combining basis is indeed unique. Finally, the design of the recombined basis is made foolproof by the procedure introduced in \cref{alg:bcs} at a very small extra cost.

Our conclusion is therefore that one should use recombined basis whenever possible to gain the substantial speed boost, particularly when a large number of differential equations are to be solved. Such scenarios include solving time-dependent PDEs \cite{che1}, solving nonlinear ODEs \cite{qin1}, computing pseudospectra \cite{den1}, etc. 

The \textsc{Julia} and \textsc{Mathematica} code used in this paper is available from \cite{qin2}.

\appendix                   
\section{Mathematica code for basis recombination}\label{sec:mcode}
The following \textsc{Mathematica} program determines the combination stencil for Chebyshev polynomials so that the new basis functions satisfy the boundary conditions at the endpoints of $[-1, 1]$. It takes as the input two lists of the orders of the boundary conditions specified by an ODE boundary value problem, one for the left boundary point and the other for the right. As the output, the code returns the expressions for the combination coefficients $\{\gamma_j^k\}_{j = 0}^{N-1}$ in the stencil matrix \cref{stencil}. In this implementation, $\gamma_N^k$ is uniformly set to $1$.

\begin{verbatim}

ClearAll["Global`*"]
(* Lists of the orders of the boundary conditions. *)
lbc = {0, 1, 2};
rbc = {1, 3, 4};

(* The order of the ODE. *)
M = Length[lbc] + Length[rbc];

(* Allocate for the stencil matrix gamma. *)
gamma = Array[\[Gamma], M, 0];

(* Initialize an empty list for the equations. *)
eqs = {};

(* Expression for the values of Chebyshev polynomials at +-1 
   (up to a scaling factor). *)
bc[deg_, ord_] := Fold[(#1*(deg^2 - (#2)^2)) &, 1, Range[0, ord - 1]]

(* Loop each boundary condition to set up equations. *)
(* Left boundary conditions. *)
For[i = 1, i <= Length[lbc], i++,
leftsign = 1;
eql = leftsign*bc[k + M, lbc[[i]]]; (* gamma[[N]] is set 1. *)
For[p = M, p >= 1, p--,
  leftsign = -leftsign; (* alternating sign *)
  eql = eql + leftsign*bc[k + p - 1, lbc[[i]]]*gamma[[p]]];
AppendTo[eqs, eql == 0];]

(* Right boundary conditions. *)
For[i = 1, i <= Length[rbc], i++,
eqr = bc[k + M, rbc[[i]]]; (* gamma[[N]] is set 1. *)
For[q = M, q >= 1, q--,
  eqr = eqr + bc[k + q - 1, rbc[[i]]]*gamma[[q]]];
AppendTo[eqs, eqr == 0];]

(* Solve the system. *)
sol = Solve[eqs[[1 ;; M]], gamma];

(* Simplify the expressions and print. *)
For[j = 1, j <= Length[sol[[1]]], j++,
  Print[Subsuperscript[\[Gamma], j-1, "k"] -> Factor[sol[[1, j, 2]]]]]
Print[Subsuperscript[\[Gamma], M, "k"] -> 1];
\end{verbatim}

\bibliographystyle{siam}
\bibliography{references.bib}

\end{document}